\nonstopmode
\documentclass[a4paper, 10pt]{amsart}
\usepackage{latexsym}
\usepackage{fancyhdr}
\usepackage{amsmath, amssymb}
\usepackage[ansinew]{inputenc}
\usepackage[all]{xy}
\usepackage{verbatim}
\usepackage{hyperref}
\swapnumbers
\theoremstyle{plain}
\newtheorem*{lemma*}{Lemma}
\newtheorem{lemma}[subsection]{Lemma}
\newtheorem*{theorem*}{Theorem}
\newtheorem{theorem}[subsection]{Theorem}
\newtheorem*{proposition*}{Proposition}

\newtheorem*{corollary*}{Corollary}
\newtheorem{corollary}[subsection]{Corollary}
\theoremstyle{definition}
\newtheorem*{definition*}{Definition}

\newtheorem*{example*}{Example}

\newtheorem*{algorithm*}{Algorithm}
\newtheorem*{remark*}{Remark}
\newtheorem{remark}[subsection]{Remark}

\newenvironment{demo}[1]{\par\smallskip\noindent{\bf #1.}}{\par\smallskip}
\numberwithin{equation}{subsection}

\def\al{\alpha}

\def\de{\delta}
\def\ep{\epsilon}

\def\rh{\rho}

\def\si{\sigma}

\def\ta{\tau}

\def\vh{\varphi}

\def\ps{\psi}
\def\om{\omega}

\def\Si{\Sigma}

\def\C{\mathbb{C}}

\def\N{\mathbb{N}}

\def\R{\mathbb{R}}

\def\cI{\mathcal{I}}
\def\cJ{\mathcal{J}}

\def\p{\partial}

\renewcommand{\Re}{\mathrm{Re}}
\renewcommand{\Im}{\mathrm{Im}}
\renewcommand{\o}{\circ}
\def\<{\langle}
\def\>{\rangle}
\def\cq{{/\!\!/}}

\let\on=\operatorname

\title[Puiseux's theorem and lifting curves]
{A generalization of Puiseux's theorem and \\ lifting curves over invariants}

\author[M.\ Losik, P.\ W.\ Michor, A.\ Rainer]
{Mark Losik, Peter W.\ Michor, and Armin Rainer}

\address{Mark Losik: Saratov State University, 
Astrakhanskaya, 83, 410026 Saratov, Russia}

\email{losikMV@info.sgu.ru}

\address{Peter W. Michor:
Fakult\"at f\"ur Mathematik, Universit\"at Wien,
Nordbergstrasse 15, A-1090 Wien, Austria} 

\email{Peter.Michor@univie.ac.at}

\address{Armin Rainer: Fakult\"at f\"ur Mathematik, Universit\"at Wien,
Nordbergstrasse 15, A-1090 Wien, Austria}

\email{armin.rainer@univie.ac.at}

\begin{document}

\begin{abstract} 
Let $\rh : G \to \on{GL}(V)$ be a rational representation of a reductive linear algebraic group $G$ defined over $\C$ on a 
finite dimensional complex vector space $V$.
We show that, for any generic smooth (resp.\ $C^M$) curve $c : \R \to V \cq G$ in the categorical quotient $V \cq G$
(viewed as affine variety in some $\C^n$)
and for any $t_0 \in \R$, there exists a positive integer $N$ such that $t \mapsto c(t_0\pm (t-t_0)^N)$
allows a smooth (resp.\ $C^M$) lift to the representation space near $t_0$. 
($C^M$ denotes the Denjoy--Carleman class associated with $M=(M_k)$, 
which is always assumed to be logarithmically convex and derivation closed). 
As an application we prove that any 
generic smooth curve in $V \cq G$ admits locally absolutely continuous (not better!) lifts.  
Assume that $G$ is finite.
We characterize curves admitting differentiable lifts.
We show that any germ of a $C^\infty$ curve which represents a lift of a germ of a quasianalytic $C^M$ curve in $V \cq G$ is 
actually $C^M$. There are applications to polar representations.
\end{abstract}

\thanks{PM was supported by the FWF-grant P21030,
AR by the FWF-grants J2771 and P22218}
\keywords{Puiseux's theorem, reductive group representations, 
invariants, regular lifting, ultradifferentiable, Denjoy--Carleman}
\subjclass[2000]{14L24, 14L30, 20G20}
\date{January 28, 2011}

\maketitle

\section{Introduction}

Let $G$ be a reductive linear algebraic group defined over $\C$ and let $\rh : G \to \on{GL}(V)$ be 
a rational representation on a finite dimensional complex vector space $V$. 
The algebra $\C[V]^G$ of invariant polynomials is finitely generated.
Choose generators $\si_1,\ldots,\si_n$, and consider the mapping $\si=(\si_1,\ldots,\si_n) : V \to \C^n$. 
We may identify $\si : V \to \si(V)$ with the morphism $\pi : V \to V \cq G$, 
where $V \cq G$ is the categorical quotient and $\pi$ is defined by the embedding $\C[V]^G \to \C[V]$.
The mapping $\si$ induces a bijection between $\si(V)$ and the closed orbits in $V$.

Let $c : \R \to V \cq G = \si(V) \subseteq \C^n$ be a smooth curve, i.e., smooth as a curve in $\C^n$.
A curve $\bar c : \R \to V$ is called a lift of $c$ to $V$, if $c = \si \circ \bar c$ and if the orbit $G.\bar c(t)$ is closed for each $t$. 
It is natural to ask whether $c$ admits lifts $\bar c$ which are regular (of some kind). 
This question is independent of the choice of generators of $\C[V]^G$, 
as any two of them differ just by a polynomial diffeomorphism.

Suppose that $G=\on{S}_n$ is acting on $V=\C^n$ by permuting the coordinates. 
The elementary symmetric functions 
$\si_j(z) = \sum_{i_1<\cdots<i_j} z_{i_1} \cdots z_{i_j}$, for $1 \le j \le n$,
generate $\C[\C^n]^{\on{S}_n}$.
By Vieta's formulas, a curve in $\C^n \cq \on{S}_n = \si(\C^n) = \C^n$ can be interpreted as a curve of monic univariate polynomials
$P(t)(z) = z^n + \sum_{j=1}^n (-1)^j a_j(t) z^{n-j}$, and, in this picture, a lift over $\si$ 
represents a parameterization of the roots of $P$. 
This special case is studied in \cite{RainerAC} and \cite{RainerQA} (see also references therein). 
It is proved that a generic smooth curve of polynomials $P$ admits a locally absolutely continuous parameterization of its roots.
Actually, any continuous choice of roots (which always exists) is locally absolutely continuous.
By \emph{generic} we mean that no two of the continuously chosen roots meet of infinite order of flatness. 
Simple examples show that one cannot expect more than absolute continuity.
This result follows from the proposition that for any $t_0$ there exists a positive integer $N$ such that 
$t \mapsto P(t_0 \pm (t-t_0)^N)$ admits smooth parameterizations of its roots near $t_0$. 

The main result of the present paper (theorem \ref{main}) is a generalization of this Puiseux type 
theorem to the framework of 
representations $\rh : G \to \on{GL}(V)$. 
We show that, for any generic smooth (resp.\ $C^M$) curve $c : \R \to V \cq G = \si(V) \subseteq \C^n$ and 
for any $t_0 \in \R$, there exists an integer $N$ such that $t \mapsto c(t_0\pm (t-t_0)^N)$
allows a smooth (resp.\ $C^M$) lift to $V$ (where $C^M$ denotes the Denjoy--Carleman class associated to $M=(M_k)$, cf.\ \ref{DC}). 
Similarly, for any holomorphic $c : \C \to V \cq G$ and any 
$z_0 \in \C$, there exists an integer $N$ such that $z \mapsto c(z_0+ (z-z_0)^N)$
allows a holomorphic lift.  
The genericity condition for $c$ is normal nonflatness \ref{gen} which, 
roughly speaking, means that $c$ does not meet lower dimensional isotropy type strata 
of $V \cq G$ with infinite order of flatness. 
It is automatically satisfied if $c$ belongs to any quasianalytic class $C^M$. 
In section \ref{secglobal} we show a global formulation of theorem \ref{main}.

As a first consequence of the main result we prove in section \ref{secAC} that any 
generic smooth curve $c : \R \to V \cq G= \si(V) \subseteq \C^n$ 
admits a locally absolutely continuous lift. 
If $G$ is finite, then any continuous lift of $c$ is locally absolutely continuous.
In general we cannot expect a better regularity than absolute continuity.
However, it is an open question whether non-generic smooth curves $c$ have locally absolutely continuous lifts.

In section \ref{secdiff} we characterize, for finite $G$, those curves in $V \cq G$ which admit differentiable lifts to $V$. 
Roughly speaking, curves must meet lower dimensional isotropy type strata of $V \cq G$ with 
sufficient order of flatness. 

We show in section \ref{secsmooth} that,
for finite $G$, any germ of a $C^\infty$ curve which represents a lift of a germ of 
a quasianalytic $C^M$ curve in $V \cq G$ is actually $C^M$. 

In section \ref{secpol} we give applications to polar representations.

In this paper we concentrate on path lifting. The problem of lifting mappings in several variables was studied in 
\cite{KLMR08} and lately in \cite{RainerRG}.

We thank the referee for pointing out that \cite{Pereira10} has recently used 
holomorphic liftings of this type to prove the Nash conjecture for quotient surface singularities.

\section{Preliminaries}

\subsection{The setting: Representations of reductive algebraic groups}
Cf.\ \cite{PV89}.
Let $G$ be a reductive linear algebraic group defined over $\C$ and let $\rh : G \to \on{GL}(V)$ be 
a rational representation on a finite dimensional complex vector space $V$. It is well-known that the algebra $\C[V]^G$ of 
$G$-invariant polynomials on $V$ is finitely generated. 
Let $\si_1,\ldots,\si_n$ be a 
system of homogeneous generators of $\C[V]^G$ with positive degrees $d_1,\ldots,d_n$.
We may consider the \emph{categorical quotient} 
$V \cq G$, i.e., the affine algebraic variety with coordinate ring $\C[V]^G$, and the 
morphism $\pi : V \to V \cq G$ defined by the embedding $\C[V]^G \to \C[V]$.
Then we can identify $\pi$ with the mapping $\si=(\si_1,\dots\si_n) : V \to \si(V) \subseteq \C^n$ and $V \cq G$ with $\si(V)$
(which we shall do consistently).
Each fiber of $\si$ contains exactly one closed orbit. If $v \in V$ and the orbit $G.v$ is closed, 
then the isotropy group $G_v = \{g \in G : g.v=v\}$ is reductive. 

\subsection{Luna's slice theorem} \label{slice}
We state a version \cite{Schwarz80} of Luna's slice theorem \cite{Luna73}.
Recall that $U$ is a \emph{$G$-saturated} subset of $V$ if $\pi^{-1}(\pi(U))=U$ 
and that a map between smooth complex algebraic varieties is \emph{\'etale} if its differential is 
everywhere an isomorphism.

\begin{theorem*}[{\cite{Luna73}}, {\cite[5.3]{Schwarz80}}] 
Let $G.v$ be a closed orbit, $v \in V$. 
Choose a $G_v$-splitting of $V \cong T_v V$ as $T_v (G.v) \oplus N_v$ and let $\vh$ denote the mapping
\[
G \times_{G_v} N_v \to V, \quad [g,n] \mapsto g(v+n).
\]
There is an affine open $G$-saturated subset $U$ of $V$ and an affine open $G_v$-saturated neighborhood 
$S_v$ of $0$ in $N_v$ such that 
\[
\vh : G \times_{G_v} S_v \to U \quad \text{and} \quad \bar \vh : (G \times_{G_v} S_v) \cq G \to U \cq G
\]
are \'etale, where $\bar \vh$ is the map induced by $\vh$.
Moreover, $\vh$ and the natural map $G \times_{G_v} S_v \to S_v \cq G_v$ induce a $G$-isomorphism
of $G \times_{G_v} S_v$ with $U \times_{U \cq G} S_v \cq G_v$.
\end{theorem*}

\begin{corollary*}[{\cite{Luna73}}, {\cite[5.4]{Schwarz80}}] 
Choose a $G$-saturated neighborhood $\overline S_v$ of $0$ in $S_v$ (classical topology) such that the canonical 
mapping $\overline S_v \cq G_v \to \overline U \cq G$ is a complex analytic isomorphism, where 
$\overline U = \pi^{-1}(\bar \vh((G \times_{G_v} \overline S_v) \cq G))$.
Then $\overline U$ is a $G$-saturated neighborhood of $v$ and $\vh : G \times_{G_v} \overline S_v \to \overline U$ is biholomorphic.
\end{corollary*}

A \emph{slice representation} of $\rh$ is a rational representation $G_v \to \on{GL}(V/T_v (G.v))$, where $G.v$ is closed.

\subsection{Luna's stratification} \label{strat}
Cf.\ \cite{Luna73}, \cite{Schwarz80}, and \cite{PV89}.
Let $v \in V$ and let $G_v$ be the isotropy group of $G$ at $v$. 
Denote by $(G_v)$ its conjugacy class in $G$, also called an \emph{isotropy class}. 
If $(L)$ is a conjugacy class of subgroups of $G$, let $(V \cq G)_{(L)}$ denote the set of points in 
$V \cq G$ corresponding to closed orbits with isotropy group in $(L)$, 
and put $V_{(L)} := \pi^{-1}((V \cq G)_{(L)})$.
Then the collection $\{(V \cq G)_{(L)}\}$ is a finite stratification of $V \cq G$ into 
locally closed irreducible smooth algebraic subvarieties.
The isotropy classes are partially ordered, namely $(H) \le (L)$ if $H$ is conjugate to a subgroup of $L$. 
If $(V \cq G)_{(L)} \ne \emptyset$, then its Zariski closure  is equal to $\bigcup_{(M) \ge (L)} (V \cq G)_{(M)} = \pi(V^L)$,
where $V^L$ is the set of $v \in V$ fixed by $L$.
There exists a unique minimal isotropy class $(H)$ corresponding to a 
closed orbit, the \emph{principal isotropy class}. Closed orbits $G.v$ with $G_v \in (H)$ are called principal.
The subset $(V \cq G)_{(H)} \subseteq V \cq G$ is Zariski open.
If we set $V_{\<H\>} := \{v \in V : G.v ~\text{closed and}~ G_v=H\}$, then $\pi$ restricts to a 
principal $(N_G(H)/H)$-bundle $V_{\<H\>} \to (V \cq G)_{(H)}$, where $N_G(H)$ denotes the normalizer of $H$ in $G$.

\subsection{Lifting curves over invariants} \label{deflift}
Let $I \subseteq \R$ be an interval.
Let $c : I \to V \cq G = \si(V) \subseteq \C^n$ be a smooth curve, i.e., smooth as a curve in $\C^n$.
A curve $\bar c : I \to V$ is called a \emph{lift} of $c$ to $V$, if $c = \si \circ \bar c$ and if the orbit $G.\bar c(t)$ is closed for each $t \in I$. 
Lifting smooth curves over invariants is independent of the 
choice of generators of $\C[V]^G$, 
as any two of them differ just by a polynomial diffeomorphism (cf.\ \cite[2.2]{KLMR06}).

\subsection{Normal nonflatness} \label{gen}
Let $s \in \N$. 
Denote by $A_s$ the union of all strata $(V \cq G)_{(L)}$ of $V \cq G$ with $\dim (V \cq G)_{(L)} \le s$, 
and by $\cI_s$ the ideal of $\C[V \cq G] = \C[V]^G$ consisting of all polynomials 
vanishing on $A_{s-1}$.
Let $I \subseteq \R$ be an open interval.
Let $c : I \to V \cq G = \si(V) \subseteq \C^n$ be
a smooth curve, $t \in I$, and $s = s(c,t)$ a minimal integer 
such that, for a neighborhood $J$ of $t$ in $I$, 
we have $c(J) \subseteq A_s$. 
The curve $c$ is called \emph{normally nonflat at $t$}, 
if there is an $f \in \cI_s$ such that $f \circ c$ is nonflat at $t$,
i.e., the Taylor series of $f \circ c$ at $t$ is not identically zero.
We say that $c$ is \emph{normally nonflat}, 
if $c$ is normally nonflat at each $t$.

Normal nonflatness is a condition which depends on the representation $\rh$. 
We shall also deal with a condition which guarantees normal nonflatness, but is independent of $\rh$.
That condition is that $c : I \to V \cq G=\si (V) \subseteq \C^n$ belongs to a quasianalytic Denjoy--Carleman class 
(cf.\ \ref{DC}).

\subsection{Denjoy--Carleman classes} \label{DC}
See \cite{Thilliez08} or \cite{KMRc} and references therein.
Let $U \subseteq \R^q$ be open. 
Let $M=(M_k)_{k \in \N}$ be a non-decreasing sequence of real numbers with $M_0=1$.
Denote by $C^M(U)$ the set of all $f \in C^\infty(U)$ such that for every compact $K \subseteq U$ there are constants $C,\rh >0$ with 
\begin{equation} \label{CM}
|\p^\al f(x)| \le C \rh^{|\al|} |\al|!\, M_{|\al|} \quad \text{ for all } \al \in \N^q \text{ and } x \in K. 
\end{equation} 
We call $C^M(U)$ a \emph{Denjoy--Carleman class} of functions on $U$.
If $M_k=1$, for all $k$, then $C^M(U)$ coincides with the ring $C^\om(U)$ of real analytic functions
on $U$. In general, $C^\om(U) \subseteq C^M(U) \subseteq C^\infty(U)$. 

We assume that $M=(M_k)$ is \emph{logarithmically convex}, i.e.,
\begin{equation} \label{logconvex}
M_k^2 \le M_{k-1} \, M_{k+1} \quad \text{ for all } k.
\end{equation}
Hypothesis \eqref{logconvex} implies that $C^M(U)$ is a ring, for all open
subsets $U \subseteq \R^q$. 
Note that definition \eqref{CM} makes sense also for mappings
$U\to \mathbb R^p$.
Then \eqref{logconvex} guarantees stability under composition and taking the inverse (i.e., the inverse function theorem holds in $C^M$).

We shall also require that $C^M$ is stable under derivation, which is
equivalent to 
\begin{equation} \label{der}
\sup_{k \in \N_{>0}} \Big(\frac{M_{k+1}}{M_k}\Big)^{\frac{1}{k}} < \infty.
\end{equation}
By the standard integral formula, stability under derivation implies that $C^M$ is closed 
under division by a coordinate (i.e., if $f \in C^M(U)$ vanishes on $\{x_i=a_i\}$, then $f(x)=(x_i-a_i) g(x)$, where $g \in C^M(U)$).

We say that $C^M$ is \emph{quasianalytic} if, 
for connected open subsets $U \subseteq \R^q$, the Taylor series homomorphism $C^M(U) \to \R[[x]], f \mapsto \hat f_a(x)$
is injective for any $a \in U$.  
Suppose that $M$ is logarithmically convex. Then,
by the Denjoy--Carleman theorem,  
$C^M$ is quasianalytic if and only if 
\begin{equation} \label{eqqa}
\sum_{k=1}^\infty \frac1{(k!M_k)^{1/k}} = \infty \quad \text{or, equivalently, } \quad \sum_{k=0}^\infty \frac{M_k}{(k+1)M_{k+1}}=\infty. 
\end{equation}

\begin{definition*}
By a {\it DC-weight sequence} we mean a non-decreasing sequence 
$M=(M_k)_{k \in \N}$ of positive numbers with $M_0=1$ 
which satisfies
\eqref{logconvex} and \eqref{der}. We say that a DC-weight sequence $M$ is \emph{quasianalytic} if \eqref{eqqa} is fulfilled.
\end{definition*}

From now on $M$ will always denote a DC-weight sequence.

Let $U \subseteq \R^q$ be open and let $f : U \to \C$ be a complex valued function. 
Then we say that $f \in C^M(U,\C)$ if $(\Re f, \Im f) \in C^M(U,\R^2)$.
If $M$ is quasianalytic, then clearly each $f \in C^M(U,\C)$ satisfying $\p^\al f(a)=0$ for all $\al \in \N^q$ is identically $0$ near $a \in U$.

Let $I \subseteq \R$ be an open interval.
A curve $c =(c_1,\ldots,c_n): I \to \C^n$ belongs to $C^M(I,\C^n)$, if each coordinate function $c_i \in C^M(I,\C)$.
Since $C^M$ is stable under composition, we may consider $C^M$ curves $c \in C^M(I,V)$ 
in vector spaces $V$. 

\subsection{Absolutely continuous functions} \label{ac}
Let $I \subseteq \R$ be an interval. A function $f : I \to \C$ is called 
\emph{absolutely continuous}, or $f \in AC(I)$, if for all $\ep > 0$ there exists a $\de >0$ 
such that $\sum_{i=1}^N (b_i-a_i) < \de$ implies $\sum_{i=1}^N |f(b_i)-f(a_i)| < \ep$,
for all sequences of pairwise disjoint subintervals $(a_i,b_i) \subseteq I$, $1 \le i \le N$.
By the fundamental theorem of calculus for the Lebesgue integral,
$f \in AC([a,b])$ if and only if there is a function $g \in L^1([a,b])$ such that 
\[
f(t) = f(a) + \int_a^t g(s) ds \quad \text{ for all } t \in [a,b]. 
\]
Then $f'=g$ almost everywhere.
Every Lipschitz function is absolutely continuous.

Gluing finitely many absolutely continuous functions provides an 
absolutely continuous function: Let $f_1 \in AC([a,b])$, $f_2 \in AC([b,c])$, and $f_1(b)=f_2(b)$.
Then $f : [a,c] \to \C$, defined by $f(t) = f_1(t)$ if $t \in [a,b]$ and $f(t)=f_2(t)$ if $t \in [b,c]$, belongs to $AC([a,c])$.
Similarly for more than two functions.

Let $\vh : I \to I$ be bijective, strictly monotone, and Lipschitz continuous. If $f \in AC(I)$ then also $f \circ \vh \in AC(I)$. 
Furthermore (cf.\ \cite[Lemma 2.5]{RainerAC}):

\begin{lemma*} 
Let $r >0$ and $n \in \N_{>0}$.
Let $f \in AC([0,r])$ (resp.\ $f \in AC([-r,0])$) and set $h(t) = f(\sqrt[n]{t})$
(resp.\ $h(t)=f(-\sqrt[n]{|t|})$). Then $h \in AC([0,r^n])$ (resp.\ $h \in AC([-r^n,0])$). 
\end{lemma*}

\begin{demo}{Proof}
There exists a function $g \in L^1([0,r])$ such that 
\[
f(t) = f(0) + \int_0^t g(s) ds
\]
for all $t \in [0,r]$. The function $(0,r^n] \to (0,r], t \mapsto \sqrt[n]{t}$, is smooth and bijective, so 
\[
\int_0^{r^n} |g(\sqrt[n]{s})| (\sqrt[n]{s})' ds = \int_0^r |g(s)| ds 
\]
and $t \mapsto g(\sqrt[n]{t}) (\sqrt[n]{t})'$ belongs to $L^1([0,r^n])$. Thus $h(t)=f(\sqrt[n]{t})$ is 
in $AC([0,r^n])$.

For the second statement consider the absolutely continuous function $f \o S|_{[0,r]}$, where $S : \R \to \R, t \mapsto -t$. 
By the above, $h_S(t) = (f \o S|_{[0,r]})(\sqrt[n]{t})$ is in $AC([0,r^n])$, and so 
$h(t) = h_S(S^{-1}|_{[-r^n,0]}(t)) = f(-\sqrt[n]{-t}) = f(-\sqrt[n]{|t|})$
is in $AC([-r^n,0])$. 
\qed\end{demo}

A curve $c =(c_1,\ldots,c_n): I \to \C^n$ is absolutely continuous, or $c \in AC(I,\C^n)$, if each coordinate function $c_i \in AC(I)$.
Since absolute continuity is invariant under linear coordinate changes, we may consider absolutely continuous curves $c \in AC(I,V)$ 
in vector spaces $V$. 

\section{An analogue of Puiseux's theorem} \label{secmain}

\subsection{Removing fixed points} \label{fix}
Let $V^G$ be the subspace of $G$-invariant vectors, and let $V'$ be a $G$-invariant complementary subspace in $V$. 
Then $V=V^G \oplus V'$, $\C[V]^G = \C[V^G] \otimes \C[V']^G$, and $V \cq G = V^G \times V' \cq G$. 
The following lemma is obvious.

\begin{lemma*}
Any smooth (resp.\ $C^M$) lift of a smooth (resp.\ $C^M$) curve $c=(c_0,c_1)$ in $V^G \times V' \cq G \subseteq \C^n$ has the form 
$\bar c = (c_0,\bar c_1)$, where $\bar c_1$ is a smooth (resp.\ $C^M$) lift of $c_1$ to $V'$. \qed
\end{lemma*}

If $c=(c_0,c_1)$ is a normally nonflat curve in $V^G \times V' \cq G$, then 
$c_1$ is a normally nonflat curve in $V' \cq G$, 
since $(V^G \times V' \cq G)_{(L)}=V^G \times (V' \cq G)_{(L)}$ (cf.\ \cite[Proposition 3.5]{AKLM00}).

\subsection{Multiplicity} \label{mult}
For a continuous real or complex valued function $f$ defined near $0$ in $\R$, 
let the \emph{multiplicity} (or \emph{order of flatness}) $m(f)$ at $0$ be the supremum 
of all integers $p$ such that $f(t)=t^p g(t)$ near $0$ for a continuous function $g$. 
Note that, if $f$ is of class $C^n$ and $m(f) < n$, then $f(t) = t^{m(f)} g(t)$ 
near $0$, where now $g$ is $C^{n-m(f)}$ and $g(0) \ne 0$. 
Similarly, one can define the multiplicity of a function at any $t \in \R$.

\begin{theorem} \label{main}
Assume that $M$ is a DC-weight sequence.
Let $I \subseteq \R$ be an open interval and let $t_0 \in I$.
Let $c = (c_1,\ldots,c_n) : I \to V \cq G = \si(V) \subseteq \C^n$ be smooth (resp.\ $C^M$). 
Suppose that $c$ is normally nonflat at $t_0$. 
Then there exists a positive integer $N$ such that $t \mapsto c(t_0 \pm (t-t_0)^N)$ allows a smooth (resp. $C^M$) lift $\bar c^\pm$ to $V$, 
locally near $t_0$.  
\end{theorem}

Note that $c$ is automatically normally nonflat at $t_0$ if $M$ is quasianalytic.

\begin{demo}{Proof}
We may assume without loss that $t_0=0$. 
Let $v \in \si^{-1}(c(0))$ such that $G.v$ is a closed orbit.
We show that there exists a positive integer $N$ such that $t \mapsto c(\pm t^N)$ has a smooth (resp.\ $C^M$) local lift $\bar c^\pm$ 
with $\bar c^\pm(0)=v$. 
Let us use the following:

\begin{algorithm*}
$(1)$ If $c(0) \in (V \cq G)_{(H)}$ is principal, then, for each $v \in V_{\<H\>} \cap \si^{-1}(c(0))$, 
a smooth (resp.\ $C^M$) lift $\bar c^{\pm}$ of $t \mapsto c(\pm t)$ to $V_{\<H\>}$ 
with $\bar c^{\pm}(0)=v$ exists, 
locally near $0$, 
since $V_{\<H\>} \to (V \cq G)_{(H)}$ is a principal $(N_G(H)/H)$-bundle (see \ref{strat}).

$(2)$ If $V^G \ne \{0\}$, we remove fixed points, by lemma \ref{fix}.

$(3)$ If $V^G = \{0\}$ and $c(0) \ne 0$ is not principal, we consider the slice representation $G_v \to \on{GL}(N_v)$.
By Luna's slice theorem \ref{slice}, the lifting problem reduces to the group $G_v$ acting on $N_v$.
As $V^G = \{0\}$, $G_v$ is a proper subgroup of $G$.
Closed $G_v$-orbits in $N_v$ correspond to closed $G$-orbits in $V$.
The stratification of $V \cq G$ in a neighborhood of $c(0)$ is naturally isomorphic
to the stratification of $N_v \cq G_v$ in a neighborhood of $0$. Since the notion of normal nonflatness is 
local, the reduced curve is normally nonflat at $0$ as well.

If $N_v^{G_v} \ne \{0\}$ we continue in $(2)$, otherwise in $(4)$.

$(4)$ If $V^G = \{0\}$ and $c(0) = 0$, we have $m(c_k) \ge 1$ for all $1 \le k \le n$. 
If $c = 0$, we lift it by $\bar c=0$ and are done. Otherwise,
put
\[
m:= \min \Big\{\frac{m(c_k)}{d_k} : 1 \le k \le n\Big\}.
\]
(Recall $d_k=\deg \si_k$.) 
Then $m$ is well-defined and finite, since $c$ is normally nonflat at $0$. 
Let $d$ be the minimal integer such that $d m \ge 1$.
Then, the multiplicity of $t \mapsto c_k(\pm t^d)$ (for $1 \le k \le n$) satisfies
\[
m(c_k(\pm t^d)) = d m(c_k) \ge d m d_k \ge d_k.
\] 
Hence, near $0$, we have $c_k(\pm t^d) = t^{d_k} c_{(1),k}^\pm(t)$, where $c_{(1),k}^\pm$ is smooth (resp.\ $C^M$), for all $k$.
Consider the smooth (resp.\ $C^M$) curve
\[
c_{(1)}^\pm(t) := (c_{(1),1}^\pm(t),\ldots,c_{(1),n}^\pm(t))
\]
which lies in $V \cq G = \si(V)$, by the homogeneity of the generators $\si_k$.
It is easy to see that $c_{(1)}^\pm$ is normally nonflat at $0$ (cf.\ \cite[Proposition 3.5]{AKLM00}).
If $t \to c_{(1)}^\pm(t)$ admits a smooth (resp. $C^M$) lift $t \mapsto \hat c_{(1)}^\pm(t)$, 
then $t \mapsto t \hat c_{(1)}^\pm(t)$ is a smooth (resp.\ $C^M$) lift of $t \mapsto c(\pm t^d)$.

Note that $m(c_{(1),k}^\pm) = d m(c_k)-d_k$ (for $1 \le k \le n$), and, thus,
\begin{equation} \label{tm}
m_{(1)} := \min  \Big\{\frac{m(c_{(1),k}^\pm)}{d_k} : 1 \le k \le n\Big\} = dm-1 < m, 
\end{equation}
by the minimality of $d$.

If $m_{(1)} = 0$ there exists some $k$ such that $c_{(1),k}^\pm(0) \ne 0$, 
and we feed $c_{(1)}^\pm$ into step $(1)$ or $(3)$.
Otherwise we feed $c_{(1)}^\pm$ into step $(4)$. 
\end{algorithm*}

Each of the steps $(1)-(3)$ either provides a required lift or reduces the lifting problem to a `smaller' space or group.
Since $m_{(1)}$ is of the form $p/d_k$, for some $1 \le k \le n$ and $p \in \N$, and by \eqref{tm}, 
also step $(4)$ is visited only finitely many times. 
So the algorithm stops and it provides a positive integer $N$ and a smooth (resp.\ $C^M$) 
lift $\bar c^\pm$ of $t \mapsto c(\pm t^N)$ near $0$  with $\bar c^\pm(0)=v$. 
\qed\end{demo}

Carrying out the obvious modifications in the proof of theorem \ref{main} we obtain:

\begin{theorem} \label{hol}
Let $U \subseteq \C$ be open and connected, and let $z_0 \in U$. 
Assume that $c : U \to V \cq G = \si(V) \subseteq \C^n$ is holomorphic.
Then there exists a positive integer $N$ such that $z \mapsto c(z_0 + (z-z_0)^N)$ allows a holomorphic lift $\bar c$ into $V$, 
locally near $z_0$. 
\qed
\end{theorem}

\begin{remark}
Note that theorem \ref{hol} is a generalization of Puiseux's theorem \cite{Puiseux1850} (see also \cite{BM90} for a modern account).
\end{remark}

\section{Smooth lifts} \label{secglobal}

\subsection{Exceptional points} \label{ex}
Let $I \subseteq \R$ be an interval.
For a curve $c : I \to V \cq G = \si(V) \subseteq \C^n$
let $E(c)$ denote the set of \emph{exceptional points}, i.e., 
the set of all $t \in I$ such that $f(c(t))=0$ for all $f \in \cI_s$  
(where $s=s(c,t)$ is defined in \ref{gen}). 

\begin{lemma*} 
If $c$ is normally nonflat, then $E(c)$ is discrete.
\end{lemma*}

\begin{demo}{Proof}
Let $t \in E(c)$. Let $s$ be minimal such that $c(J) \subseteq A_s$ for a neighborhood $J$ of $t$. 
By assumption, there exists a $f \in \cI_s$ such that $f \circ c$ is nonflat at $t$. Since $t \in E(c)$, 
we have $f(c(t))=0$. The mean value theorem implies that $t$ is an isolated zero of $f \circ c$ 
(see e.g.\ \cite[Lemma 3.2]{RainerAC} for details).
So $t$ is isolated in $E(c)$. 
\qed\end{demo}

\begin{lemma} \label{unique}
Let $c : I \to V \cq G = \si(V) \subseteq \C^n$ be  a smooth (resp.\ $C^M$) curve which is normally nonflat at $t_0 \in I$.
Suppose that $\bar c_1, \bar c_2 : J \to V$ are smooth (resp.\ $C^M$) lifts of $c$ on an open subinterval $J \subseteq I$ containing $t_0$.
Then there exists a smooth (resp.\ $C^M$) curve $g$ in $G$ defined near $t_0$ such that $\bar c_1(t) = g(t).\bar c_2(t)$ for all $t$ near $t_0$.
\end{lemma}

\begin{demo}{Proof}
The proof is essentially the same as in \cite[3.8]{AKLM00}.

Without loss $0 \in I$, $t_0=0$, and $\bar c_1(0)=\bar c_2(0)=:v$.
Consider the projection 
$p : G.S_{v} \cong G \times_{G_{v}} S_{v} \to G/G_{v}$
of the fiber bundle associated to the principal bundle $G \to G/G_v$. Then, for $t$ near $0$, 
$p \o \bar c_1$ and $p \o \bar c_2$ admit smooth (resp.\ $C^M$ lifts) $g_1$ and $g_2$ to $G$ with $g_1(0)=g_2(0)=e$, 
and $t \mapsto g_j^{-1}(t).\bar c_j(t)$ form smooth (resp.\ $C^M$) lifts of $c$ lying in $S_v$. 
This reduces the problem to the group $G_v$.
If $c(0)$ is principal, then $G_v$ acts trivially on $N_v$, and the two lifts coincide.
Otherwise, we may remove fixed points (by \ref{fix}) and assume $c(0)=0$ and $V^G = \{0\}$.

If $c=0$ identically, the statement is trivial. Otherwise, $m(\bar c_1) = m(\bar c_2) = r < \infty$, since $c$ is normally nonflat at $0$.
(Here we use the obvious notion of multiplicity for a curve in a vector space.) Then $t \mapsto t^{-r} \bar c_j(t)$ are smooth 
(resp.\ $C^M$) lifts of the smooth (resp.\ $C^M$) curve 
\[
c_{(r)}(t) := \si(t^{-r} \bar c_j(t)) = (t^{-d_1 r}c_1(t),\ldots,t^{-d_n r}c_n(t)).
\]
If we find $g(t) \in G$ taking $t^{-r} \bar c_1(t)$ to $t^{-r} \bar c_2(t)$, then we also have $\bar c_1(t) = g(t).\bar c_2(t)$.
Since $c_{(r)}(0) \ne 0$, we are done by induction. 
\qed\end{demo}

\subsection{Notation}
Let $I \subseteq \R$ be an open interval, $s \in I$, and $N \in \N_{>0}$. We denote by $\ps_{s,I}^{N,\pm}$ the mapping 
\begin{align*}
\ps_{s,I}^{N,\pm} : \{t \in \R : s \pm (t-s)^{N} \in I\} \to I :
t \mapsto s \pm (t-s)^{N}.
\end{align*}

\begin{theorem}
Assume that $M$ is a non-quasianalytic DC-weight sequence.
Let $I \subseteq \R$ be an open interval.
Let $c : I \to V \cq G = \si(V) \subseteq \C^n$ be smooth (resp.\ $C^M$) and normally nonflat.
For each compact subinterval $J \subseteq I$, 
there exist a neighborhood $I_1$ of $J$ and finitely many $s_1,\ldots,s_m$ and $N_1,\ldots,N_m$ such that 
$c \o \ps_{s_1,I_1}^{N_1,\pm} \o \cdots \o \ps_{s_m,I_m}^{N_m,\pm}$ admits a global smooth (resp.\ $C^M$) lift
(for each choice of signs),
where $I_{i+1} =(\ps_{s_i,I_i}^{N_i,\pm})^{-1}(I_i)$ for $i\ge 1$. 
(By construction, $I_1$ is covered by the images under $\ps_{s_1,I_1}^{N_1,\pm} \o \cdots \o \ps_{s_m,I_m}^{N_m,\pm}$ 
of the domains of $c \o \ps_{s_1,I_1}^{N_1,\pm} \o \cdots \o \ps_{s_m,I_m}^{N_m,\pm}$.)
\end{theorem}

\begin{demo}{Proof}
The algorithm in \ref{main} shows that $c$ admits smooth (resp.\ $C^M$) lifts near each $t \not\in E(c)$.

Let $J$ be fixed. By lemma \ref{ex}, the set $E(c) \cap J$ is finite; let $t_1 < t_2 < \ldots <t_m$ denote its elements.
Let $I_1$ be an open neighborhood of $J$ such that $E(c) \cap J=E(c) \cap I_1$.
By theorem \ref{main}, there is an $N_1 \in \N_{>0}$ such that $c \o \ps_{t_1,I_1}^{N_1,\pm}$ admits a smooth (resp.\ $C^M$) lift near $t_1$. 

If $N_1$ is odd, it suffices to consider $\ps_{t_1,I_1}^{N_1,+}$. 
If $N_1$ is even, then $\ps_{t_1,I_1}^{N_1,-}$ maps onto $I_1 \cap (\infty,t_1]$,  $E(c \o \ps_{t_1,I_1}^{N_1,-})=\{t_1\}$, 
and $c \o \ps_{t_1,I_1}^{N_1,-}$ admits smooth (resp.\ $C^M$)
lifts near every point in its domain.
So it is no restriction to just treat $c \o \ps_{t_1,I_1}^{N_1,+}$, which
admits smooth (resp.\ $C^M$) lifts near each $t < t_1 + (t_2-t_1)^{1/N_1}=:s_2$.
By theorem \ref{main}, there is an $N_2 \in \N_{>0}$ such that $c \o \ps_{t_1,I_1}^{N_1,+} \o \ps_{s_2,I_2}^{N_2,\pm}$ 
admits a smooth (resp.\ $C^M$) lift near $s_2$.

Repeating this process $m$ times, we obtain a family of curves of the required form 
each of which admits local smooth (resp.\ $C^M$) lifts near any point in its domain.
The local lifts can be glued to a global lift, using lemma \ref{unique} as in \cite[4.1]{AKLM00}. 
(Here we use that, if $M$ is non-quasianalytic, then $C^M$ cutoff functions exist.)
\qed\end{demo}

\section{Absolutely continuous lifts} \label{secAC}

\begin{theorem} \label{pathfG}
Let $G$ be finite.
Let $I \subseteq \R$ be an interval and let  $c : I \to V \cq G = \si(V) \subseteq \C^n$ be continuous. Then there exists a 
continuous lift $\bar c : I \to V$ of $c$.
\end{theorem}

\begin{demo}{Proof}
Without loss we may assume $V^G=\{0\}$ (cf.\ \ref{fix}).

We use induction on $|G|$. 
If $G = \{e\}$ put $\bar c := c$.
So let us assume that the theorem is proved for complex finite dimensional representations 
of proper subgroups of $G$.
 
Let $c : I \to V \cq G = \si(V) \subseteq \C^n$ 
be continuous.  
The set $c^{-1}(0)$ is closed in $I$ and, 
thus, 
$c^{-1}(\si(V) \setminus \{0\})$ is an at most countable disjoint union of open subintervals $J \subseteq I$
which are maximal with respect to not containing zeros of $c$.

Let one subinterval $J$ be fixed. 
Since $V^G = \{0\}$, for all $v \in V \setminus \{0\}$ we have $G_v \ne G$. 
By induction hypothesis and by \ref{slice}, 
we find local continuous lifts of $c$ near any 
$t \in J$ and through all $v \in \si^{-1}(c(t))$. 
Suppose $\bar c_1 : J \supseteq J_1 \rightarrow V \setminus \{0\}$ 
is a local continuous lift of $c$ with maximal domain $J_1$, 
where, say, the right endpoint $t_1$ of $J_1$ lies in $J$.
Then there exists a local continuous lift $\bar c_2$ 
of $c$ near $t_1$, and there is a $t_0 < t_1$ such that both $\bar c_1$ and $\bar c_2$ 
are defined near $t_0$. 
Since $\si(\bar c_1(t_0)) = \si(\bar c_2(t_0))$,   
there must exist a $g \in G$ such that 
$\bar c_1(t_0) = g.\bar c_2(t_0)$. 
But then $\bar c_{12}(t) := \bar c_1(t)$ for $t \le t_0$ and $\bar c_{12}(t) := g.\bar c_2(t)$ for $t \ge t_0$ 
is a continuous lift of $c$ defined on a larger 
interval than $J_1$. 
We can conclude that there exist continuous lifts of $c$ on $J$.

Let $\bar c$ be a continuous lift of $c$ on $I \setminus c^{-1}(0)$ and extend it to $I$, by putting 
$\bar c(t) := 0$ for $t \in c^{-1}(0)$. It remains to show that $\bar c$ is continuous on $c^{-1}(0)$.
Let $t_0 \in c^{-1}(0)$ and $I \ni t_n \to t_0$. 
Since $G$ is finite, each neighborhood of $0$ contains a $G$-invariant neighborhood of $0$.
So $\si(\bar c(t_n)) = c(t_n) \to c(t_0) = 0$ 
implies $\bar c(t_n) \to 0$.
\qed\end{demo}

\begin{remark} \label{rmkp}
We do not know whether paths are continuously liftable in general 
(but consider theorem \ref{aclifts}).
The proof in \ref{pathfG} does not work, since in general we cannot expect that each 
neighborhood of $0$ contains a $G$-invariant neighborhood of $0$. This is shown by the following example:
Let $G = \C^*$ act on $V=\C^2$ by $g.(x,y)=(gx,g^{-1}y)$. Then $\C[V]^G$ is generated by $\si(x,y) = xy$ and $V \cq G \cong \C$.
For each $z \in \C^*$ the fiber $\si^{-1}(z)$ represents a principal orbit. The fiber $\si^{-1}(0)$ consists of three orbits:
the closed orbit $\{0\}$, $G.(1,0)=\C^* \times \{0\}$, and $G.(0,1)= \{0\} \times \C^*$. 
However: Each path $c$ in $V \cq G$ has a continuous lift $\bar c=(\sqrt{c},\sqrt{c})$ to $V$.
If $c$ is sufficiently smooth, then $\bar c$ is even locally absolutely continuous. 
This example is \emph{polar}. We shall see in section \ref{secpol} that for polar representations paths are liftable.
\end{remark}

\begin{lemma} \label{Lip}
Let $G$ be finite. Choose some norm $\|.\|$ on $V$.
Let $I \subseteq \R$ be an interval and let $c : I \to V \cq G = \si(V) \subseteq \C^n$ be continuous.
If there is a Lipschitz lift of $c$, then any continuous 
lift of $c$ is Lipschitz.
\end{lemma}

\begin{demo}{Proof}
Let $\bar c$ and $\tilde c$ be continuous lifts of $c$, and suppose that $\tilde c$ is Lipschitz on $I$.
Let $t<s$ be in $I$.
There is a $g_0 \in G$ such that $\bar c(t)=g_0.\tilde c(t)$. Now let $t_1$ be the 
maximum of all $r\in [t,s]$ such that $\bar c(r)=g_0.\tilde c(r)$. If $t_1<s$ then 
$\bar c(t_1)=g_1.\tilde c(t_1)$ for some $g_1 \in G \setminus \{g_0\}$. 
(Namely: For $r_n \searrow t_1$ we have $\bar c(r_n)=g_n.\tilde c(r_n)$ with $g_n \in G \setminus \{g_0\}$. 
By passing to a subsequence, $\bar c(r_n)=g_1.\tilde c(r_n)$ with $g_1 \in G \setminus \{g_0\}$, hence the assertion.)
Let $t_2$ be the 
maximum of all $r\in [t_1,s]$ such that $\bar c(r)=g_1.\tilde c(r)$. If $t_2<s$ then 
$\bar c(t_2)=g_2.\tilde c(t_2)$ for some $g_2 \in G \setminus \{g_0,g_1\}$. 
And so on until $s=t_k$ for some $k\le |G|$. Then we have (where $t_0=t$)
\begin{align*}
\frac{\|\bar c(s)-\bar c(t)\|}{s-t} 
&\le \sum_{j=0}^{k-1} \frac{\|g_j.\tilde c(t_{j+1})-g_j.\tilde c(t_j)\|}{t_{j+1}-t_j}
\cdot\frac{t_{j+1}-t_j}{s-t}
\le C \cdot M_G,
\end{align*}
where $C$ is the Lipschitz constant of $\tilde c$ and $M_G$ denotes the maximum of the 
operator norms of the linear transformations $g \in G$.
\qed\end{demo}

\begin{theorem} \label{aclifts}
Let $I \subseteq \R$ be an open interval and let $c : I \to V \cq G = \si(V) \subseteq \C^n$ be smooth and normally nonflat. 
Then:
\begin{enumerate}
\item There exists a locally absolutely continuous lift $\bar c$ of $c$, i.e., $\bar c \in  AC_{\on{loc}}(I,V)$. 
\item If $G$ is finite, any continuous lift of $c$ is locally absolutely continuous. 
\end{enumerate}
\end{theorem}

\begin{demo}{Proof}
(1) First we show: 
{\it For each $t_0 \in I$, there exist a neighborhood $J \subseteq I$ of $t_0$ and a lift $\bar c \in AC(J, V)$ of $c$ on $J$.}  
Without restriction we may assume that $0 \in I$ and $t_0=0$. 
By theorem \ref{main}, there is a positive integer $N$ and a neighborhood $J_N$ of $0$ such that 
$t \mapsto c(\pm t^N)$ allows a smooth lift $\tilde c^\pm$ on $J_N$ 
such that $\tilde c^-(0)= \tilde c^+(0)$. 
Let $J := \{t \in I : \pm \sqrt[N]{|t|} \in J_N\}$.
We may assume that $J$ is compact (by shrinking $J_N$ if necessary). 
Let $J_{\ge 0} = \{t \in J : t \ge 0\}$ and $J_{\le 0} = \{t \in J : t \le 0\}$. 
We define
\[
\bar c(t) := \left\{
\begin{array}{ll}
\tilde c^+(\sqrt[N]{t}) & ~\text{for}~ t \in J_{\ge 0} \\
\tilde c^-(-\sqrt[N]{|t|}) & ~\text{for}~ t \in J_{\le 0} ~\text{if}~ N ~\text{is even}~\\
\tilde c^+(-\sqrt[N]{|t|}) & ~\text{for}~ t \in J_{\le 0} ~\text{if}~ N ~\text{is odd}.
\end{array}
\right.
\]
Then $\bar c$ is a lift of $c$ on $J$ and $\bar c \in AC(J,V)$, by lemma \ref{ac}.

We can glue the local absolutely continuous lifts of $c$: 
Assume that $\bar c_1$ and $\bar c_2$ are local lifts of $c$ with domains of definition 
which have a common point $t_1$ and none of them is a subset of the other. 
Since both $G.\bar c_1(t_1)$ and $G.\bar c_2(t_1)$ are closed, we have $G.\bar c_1(t_1)=G.\bar c_2(t_1)$. 
Hence we can glue $\bar c_1$ and $\bar c_2$ at $t_1$ by applying to one of them a fixed transformation of $G$. 
Eventually, we obtain a global continuous lift $\bar c : I \to V$ of $c$ which is absolutely continuous on each 
compact subinterval of $I$. 

(2) Assume that $G$ is finite and let $\bar c$ be a continuous lift of $c$.
We show that each $t_0 \in I$ has a neighborhood $J$ such that $\bar c \in AC(J,V)$. 
Without restriction assume that $t_0=0$. 
By theorem \ref{main},
we find a positive integer $N$ and a neighborhood $J_N$ of $0$ such that 
$t \mapsto c(\pm t^N)$ allows a smooth lift $\tilde c^\pm$ on $J_N$. 
Another continuous lift is provided by $t \mapsto \bar c(\pm t^N)$.
By lemma \ref{Lip}, the lift $t \mapsto \bar c(\pm t^N)$ is actually Lipschitz, thus, absolutely continuous. 
Using lemma \ref{ac}, we can conclude in a similar way as in (1) that $\bar c \in AC(J,V)$.
This completes the proof.
\qed\end{demo}

\begin{remark}
The conclusion in theorem \ref{aclifts} is best possible.
In general
there is no lift with first order derivative in $L^p_{\on{loc}}$
for any $1 < p \le \infty$.
A counter example is the curve of polynomials
\[ 
P(t)(z) = z^n- t, \quad  t \in \R,
\]
which is a curve in $\C^n \cq \on{S}_n = \C^n$.
If $n \ge \frac{p}{p-1}$, for $1 < p < \infty$,
and if $n \ge 2$, for $p=\infty$, no parameterization of the roots (lift) of $P$ has first order derivative in $L^p_{\on{loc}}$.
\end{remark}

\section{Differentiable lifts} \label{secdiff}

\subsection{Notation} \label{J}
The set $A_s$ (cf.\ \ref{gen}) is Zariski closed. 
Let $\cJ_{s+1}$ be an ideal in $\C[V \cq G] = \C[V]^G$ defining $A_s$. Since the action of $G$ commutes with 
homotheties, we can assume that $\cJ_{s+1}$ is generated by homogeneous elements.

Let $d=d(\rh)$ be the maximum of all $d_i = \deg \si_i$, where $\si_1,\ldots,\si_n$ is a minimal system of homogeneous generators of $\C[V]^G$.
Let $D=D(\rh)$ be the maximum of all $d(\rh')$, where $\rh'$ is any slice representation of $\rh$. 
There are only finitely many isomorphism types of such $\rh'$. 

\begin{lemma} \label{mult1}
Let $G$ be finite. Assume that $V^G=\{0\}$.
Let $I \subseteq \R$ be an open interval containing $0$. 
Let $c=(c_1,\dots,c_n) : I \to V \cq G = \si(V) \subseteq \C^n$ be $C^d$.
Then the following conditions are equivalent:
\begin{enumerate}
\item $m(c_k) \ge d_k$, for all $1 \le k \le n$.
\item $m(f \o c) \ge \deg f$, for all homogeneous $f \in \cJ_1$.
\end{enumerate}
\end{lemma}

\begin{demo}{Proof}
Since $V^G=\{0\}$, we have $A_0=V(\cJ_1)=\{0\}$.

$(1) \Rightarrow (2)$
The continuous curve $c_{(1)}(t):=(t^{-d_1} c_1(t),\ldots,t^{-d_n} c_n(t))$ admits a 
continuous lift $\bar c_{(1)}$ (by theorem \ref{pathfG}). Then $t \mapsto t \bar c_{(1)}(t)$ is a continuous lift of $t \mapsto c(t)$, 
and $(2)$ follows from the homogeneity of $f$.

$(2) \Rightarrow (1)$
Since $\cJ_1$ is generated by homogeneous elements, $(2)$ implies $c(0)=0$.
So there exist positive integers $m_k$ and continuous $c_{m,k}$ such that
$c_k(t)=t^{m_k} c_{m,k}(t)$, for $1 \le k \le n$. Assume, for contradiction, that $m_k=m(c_k) < d_k$ for some $k$. Then
\[
m:=\min \Big\{ \frac{m_1}{d_1},\ldots,\frac{m_n}{d_n} \Big\} < 1.
\]
For $t \ge 0$, consider 
\[
c_{(m)}(t):=(t^{m_1-d_1 m} c_{m,1}(t),\ldots,t^{m_n-d_n m} c_{m,n}(t)),
\]
which is a continuous curve in $V \cq G=\si(V)$, and hence has a continuous lift $\bar c_{(m)}$ (by theorem \ref{pathfG}).
Then $t \mapsto t^m \bar c_{(m)}(t)$ is a continuous lift of $t \mapsto c(t)$.
Let $f \in \cJ_1$ be homogeneous. By (2),
\[
m(f \o c_{(m)}) = m(f \o c) - m \deg f \ge (1-m) \deg f > 0.
\]
So $c_{(m)}(0)=0$ (by the same arguments as above). But that is a contradiction for those $k$ with $m(c_k)=m_k=d_k m$.
\qed\end{demo}

\subsection{$1$-flatness} \label{1flatness}
Let $G$ be finite. 
For $s \in \N$, put $\bar A_s := \pi^{-1}(A_s)$. 
Then $\cJ_{s+1} \subseteq \C[V]^G$ is a defining ideal of $\bar A_s$.
Let $I \subseteq \R$ be an open interval.
Let $c : I \to V \cq G = \si(V) \subseteq \C^n$ be continuous and let $\bar c : I \to V$ be a continuous lift of $c$. 
Let $t \in I$ and let $r=r(c,t)$ be the minimal integer such that $c(t) \in A_r$.
The lift $\bar c$ is called \emph{$1$-flat at $t$} 
if the multiplicity of $f \circ \bar c$ at $t$ is $\ge \deg f$ for all homogeneous $f \in \cJ_{r+1}$.

Assume that $0 \in I$ and $c : I \to V \cq G = \si(V) \subseteq \C^n$ is $C^D$.
Choose $v \in \si^{-1}(c(0))$.  
Let $\ta_1,\ldots,\ta_m$ be a minimal system of homogeneous generators of $\C[N_v]^{G_v}$ 
with degrees $d'_1,\ldots,d'_m$, 
and consider $\ta=(\ta_1,\ldots,\ta_m) : N_v \to \C^m$. 

\begin{lemma*} 
The following conditions are equivalent:
\begin{enumerate}
\item There is a continuous lift $\bar c$ with $\bar c(0)=v$ and differentiable at $0$.
\item There is a continuous lift $\bar c$ with $\bar c(0)=v$ and $1$-flat at $0$. 
\item The curve $c$, considered as curve in $\ta(N_v)$ with $c(0)=0$, is of the form 
$c(t) = (t^{d'_1} c_{(1),1}(t),\ldots,t^{d'_m} c_{(1),m}(t))$ near $0$, where $c_{(1),k}$, for $1 \le k \le m$, is continuous.
\end{enumerate}
\end{lemma*}

\begin{demo}{Proof}
$(1) \Rightarrow (2)$
Let $r=r(c,0)$ and $f \in \cJ_{r+1}$ homogeneous. 
Then $f(\bar c(0))=0$.  
By assumption, $\bar c(t) = t \tilde c(t)$ such that $\tilde c$ is continuous.
Thus $m(f \o \bar c) \ge \deg f$.

$(2) \Rightarrow (3)$
As $G$ is finite, $N_v=V$. So we may consider $\bar c$ as lift of the curve $c$ in $N_v \cq G_v = \ta(N_v)$ with $\bar c(0)=0 \in N_v$.
Then $\bar c$ is still $1$-flat at $0$ by \ref{slice}, since the notion of $1$-flatness is local. 
We may also assume that $N_v^{G_v}=\{0\}$ by \ref{fix}. So lemma \ref{mult1} implies $(3)$.

$(3) \Rightarrow (1)$
Suppose that $c(t) = (t^{d'_1} c_{(1),1}(t),\ldots,t^{d'_m} c_{(1),m}(t))$ 
near $0$ with continuous $c_{(1),k}$, for $1 \le k \le m$. 
Then $c_{(1)} = (c_{(1),1},\ldots,c_{(1),m})$ is a continuous curve in $\ta(N_v)$ which allows a continuous
lift $\bar c_{(1)}$. It follows that $\bar c(t) := t \bar c_{(1)}(t)$ is a continuous lift of $c$ near $0$ which is 
differentiable at $0$.
\qed\end{demo}

\subsection{Differentiable lifts} \label{diff1flat}
We say that $c$ is \emph{$1$-flat} at $0$ if one of the equivalent conditions in lemma \ref{1flatness} is satisfied.
Similarly for any $t \in I$. The curve $c$ is called \emph{$1$-flat}, if it is $1$-flat at each $t$.

\begin{theorem*}  
Let $G$ be finite.
Let $I \subseteq \R$ be an open interval and let $c : I \to V \cq G = \si(V) \subseteq \C^n$ be $C^D$.
Then the following conditions are equivalent:
\begin{enumerate}
\item There exists a global differentiable lift $\bar c$ of $c$.
\item $c$ is $1$-flat.
\end{enumerate}
\end{theorem*}

\begin{demo}{Proof}
$(1) \Rightarrow (2)$ is clear by lemma \ref{1flatness}. 

$(2) \Rightarrow (1)$ 
We use induction on $|G|$. 
There is nothing to prove if $G$ is trivial.
So let us assume that the statement of the theorem holds for representations of proper subgroups of $G$.

We may suppose that $V^G=\{0\}$. 
Consider the set $F := c^{-1}(0)$. 
Then $F$ is closed and its complement $I \setminus F$ is a 
countable union of open subintervals whose boundary points lie in $F$. 

Let $J$ denote one such interval. \emph{We claim that there exists a differentiable lift $\bar c$ of $c$ on $J$.} 
For each $t_0 \in J$ and each $v \in \si^{-1}(c(t_0))$, 
we have $|G_v|<|G|$, since $V^G=\{0\}$. 
By induction hypothesis and \ref{slice}, we find differentiable lifts of $c$, locally near any $t_0 \in J$. 
Let $\bar c$ be a differentiable lift of $c$ defined on a maximal 
subinterval $J' \subseteq J$. Suppose for contradiction that the right (say) endpoint $t_1$ of $J'$ 
belongs to $J$. Then there exists a differentiable lift 
$\tilde c$ of $c$, locally near $t_1$, and there is a 
$t_0 < t_1$ such that both $\bar c$ and $\tilde c$ are defined near $t_0$. 
Let $(t_m)$ be a sequence with $t_m \to t_0$. For each $m$ there exists a $g_m \in G$ such that 
$\bar c(t_m)=g_m.\tilde c(t_m)$. 
By passing to a subsequence (again denoted by $(t_m)$) 
we obtain $\bar c(t_m)=g.\tilde c(t_m)$ for a fixed $g$ and for all $m$. 
Then $\bar c(t_0) = \lim_{t_m \to t_0} \bar c(t_m) = g.(\lim_{t_m \to t_0} \tilde c(t_m)) = g.\tilde c(t_0)$ and 
\[
\bar c'(t_0) = \lim_{t_m \to t_0} \frac{\bar c(t_m)-\bar c(t_0)}{t_m-t_0} 
= \lim_{t_m \to t_0} \frac{g.\tilde c(t_m)-g.\tilde c(t_0)}{t_m-t_0}
= g.\tilde c'(t_0).
\] 
Hence, we can extend the lift $\bar c$ by $g.\tilde c$ for $t \ge t_0$ beyond $t_1$. This contradicts maximality.

Let us extend $\bar c$ to the closure of $J$, by setting it $0$ at the 
endpoints of $J$.  
Let $t_0$ denote the right (say) endpoint of $J$. 
By lemma \ref{1flatness}, there exists a local continuous 
lift $\tilde c$ of $c$ near $t_0$ which is 
differentiable at $t_0$. 
Let $(t_m)$ be a sequence with $t_m \nearrow t_0$.  
By passing to a subsequence, we may assume that 
$\bar c(t_m)=g.\tilde c(t_m)$ for a fixed $g \in G$ and for all $m$. 
Then $\lim_{t_m \nearrow t_0} \bar c(t_m) = g.(\lim_{t_m \nearrow t_0} \tilde c(t_m)) = g.\tilde c(t_0)= 0$ and 
\[
\lim_{t_m \nearrow t_0} \frac{\bar c(t_m)}{t_m-t_0} 
= \lim_{t_m \nearrow t_0} \frac{g.\tilde c(t_m)}{t_m-t_0}
= g.\tilde c'(t_0).
\] 
It follows that the set of accumulation points of $\bar c(t)/(t-t_0)$, as 
$t \nearrow t_0$, lies in the $G$-orbit through $\tilde c'(t_0)$. 
Since this orbit is finite, lemma \ref{acc} below implies 
that the limit $\lim_{t \nearrow t_0} \bar c(t)/(t-t_0)$ exists. 
Thus the one-sided derivative of $\bar c$ at $t_0$ exists. 

For isolated points in $F$ we can apply a 
fixed transformation from $G$ to one of the neighboring 
differentiable lifts in order to glue them differentiably (by the arguments in the previous paragraph). 
So we have found a differentiable lift $\bar c$ 
of $c$ defined on $I \setminus F'$, where $F'$ denotes the set of 
accumulation points of $F$.

Let us extend $\bar c$ by $0$ on $F'$. Then $\bar c$ is a global differentiable 
lift of $c$, since any lift is differentiable at points $t' \in F'$. For: 
It is clear that the derivative at $t'$ of any differentiable lift has to be $0$. 
Let $\tilde c$ be a local continuous lift which is differentiable at $t'$, provided by 
lemma \ref{1flatness}.  
As above we may conclude that 
the set of accumulation points of $\bar c(t)/(t-t')$, as 
$t \to t'$, lies in $G.\tilde c'(t') = G.0=\{0\}$. 
\qed\end{demo}

\begin{remark*}
If $\rh : G \to \on{O}(V)$ is a real finite dimensional representation of a 
compact Lie group $G$, then condition $(2)$ in theorem \ref{diff1flat} is automatically satisfied. 
Moreover, $D=d$ by \cite[Lemma 2.4]{KLMR06}.
So any $C^d$ curve $c$ in the orbit space $V/G = \si(V) \subseteq \R^n$ admits a global differentiable lift to $V$.
That case is treated in \cite{KLMR05}. 
\end{remark*}

\begin{lemma}[{\cite[Lemma 4.3]{KLMR05}}] \label{acc}
Let $X$ be a compact metric space and let $c : (a,b) \to X$ be continuous. 
Then the set of all accumulation points of $c(t)$ as $t \searrow a$ is connected.
\end{lemma}

\section{Smooth lifts of quasianalytic $C^M$-curves are $C^M$} \label{secsmooth}

We prove an analogue of \cite[Theorem 2]{Thilliez10}.

\begin{theorem} \label{smoothroots}
Let $G$ be finite.
Assume that $M$ is a quasianalytic DC-weight sequence.
Let $c : \R,0 \to V \cq G = \si(V) \subseteq \C^n$ be a germ of a $C^M$ curve.
If $\bar c : \R,0 \to V$ is a $C^\infty$ germ such that $\si \o \bar c = c$, then $\bar c$ is a $C^M$ germ.
\end{theorem}

\begin{demo}{Proof}
Choose a basis in $V$ and identify $V =\C^m$. Denote by $\on{pr}_i$ the projection to the $i$-th coordinate.
Set 
\[
P(x)(z) = \prod_{i=1}^m \prod_{g \in G} (z-\on{pr}_i(g.x)) 
= z^{m|G|}+\sum_{j=1}^{m|G|} (-1)^j a_j(x) z^{m|G|-j}.
\]
By construction the coefficients $a_j$ (for $1 \le j \le m|G|$) of $P$ belong to $\C[V]^G$.
Thus there exist polynomials $p_j$ such that $a_j = p_j \o \si$ for all $j$.
It follows that $P \o \bar c$ is a polynomial, whose coefficients are $C^M$ germs at $0 \in \R$ and which admits 
a $C^\infty$ parameterization $\on{pr}_i(g.\bar c)$ (for $1 \le i \le m$ and $g \in G$) of its roots.  
By \cite[Theorem 2]{Thilliez10}, each $\on{pr}_i(g.\bar c)$ is actually $C^M$, hence, $\bar c$ is $C^M$. 
\qed\end{demo}

\begin{remark} \label{rmk1}
If $G$ is not finite the theorem is false:
For any $C^M$ curve germ $c$ in $V \cq G$ of the example in \ref{rmkp} admitting a $C^M$-lift $\bar c$ we may 
find a $C^\infty$ curve germ $g$ in $G$ such that $g.\bar c$ is not $C^M$. 
But we may ask whether a $C^M$ curve germ admitting a $C^\infty$ lift admits also a $C^M$ lift.
\end{remark}

It seems to be unclear whether the multivariable version of \cite[Theorem 2]{Thilliez10} and thus of theorem \ref{smoothroots} 
is true. But in the real analytic case we have the following analogue of \cite[Theorem 2]{Siciak70}. 

\begin{corollary}
Let $G$ be finite.
Let $f : \R^q,0 \to V \cq G = \si(V) \subseteq \C^n$ be a germ of a $C^\om$ mapping.
If $\bar f : \R^q,0 \to V$ is a $C^\infty$ germ such that $\si \o \bar f = f$, then $\bar f$ is a $C^\om$ germ.
\end{corollary}

\begin{demo}{Proof}
If $q=1$ this is a special case of theorem \ref{smoothroots}. A smooth mapping which is $C^\om$ along $C^\om$ curves 
(affine lines suffice) is $C^\om$ (e.g.\ \cite[Theorem 1]{Siciak70}).   
\qed\end{demo}

\section{Generalization to polar representations} \label{secpol}

\subsection{Polar representations}
Cf.\ \cite{DK85}. 
Let $G$ be a reductive linear algebraic group defined over $\C$ and let $\rh : G \to \on{GL}(V)$ be 
a rational representation on a finite dimensional complex vector space $V$.
Let $v\in V$ be such that $G.v$ is closed and
consider the subspace $\Si_v =\{x \in V : \mathfrak g.x \subseteq \mathfrak g.v\}$. 
Then for each $x \in \Si_v$ the orbit $G.x$ is closed.
The representation $\rh$ is called \emph{polar} if there is a $v \in V$ with $G.v$ closed such that $\on{dim} \Si_v = \on{dim} \C[V]^G$. 
Such $\Si_v$ is called a \emph{Cartan subspace}.
Any two Cartan subspaces are conjugated. 
All closed orbits in $V$ intersect $\Si_v$. 
The \emph{generalized Weyl group} 
\[
W(\Si_v)=\{g\in G : g.\Si_v=\Si_v\}/\{g \in G : g.x=x \text{ for all } x\in \Si_v\}
\] 
is finite.
The restriction map $f\mapsto f|_{\Si_v}$ induces an isomorphism $\C[V]^G \rightarrow \C[\Si_v]^{W(\Si_v)}$. 
So we have $V \cq G = \si(V) = \si_{\Si_v}(\Si_v) = \Si_v \cq W(\Si_v)$. 

\begin{theorem}
Let $\rh : G \to \on{GL}(V)$ be polar, $\Si$ a Cartan subspace, and $W=W(\Si)$.
Let $D=D(\rh)$ be as in \ref{J}.
Let $I \subseteq \R$ be an open interval and 
consider a curve $c : I \to V \cq G = \si(V) \subseteq \C^n$. Then:
\begin{enumerate}
\item If $c$ is continuous, then there exists a continuous lift $\bar c : I \to \Si \subseteq V$ of $c$.
\item If $c$ is smooth and normally nonflat, then any continuous lift $\bar c$ of $c$ to $\Si$ is locally absolutely continuous.
\item If $c$ is $C^D$, then there exists a differentiable lift of $c$ to $\Si$ if and only if $c$ is $1$-flat with respect to $W \to \on{GL}(\Si)$. 
\item Let $M$ be a quasianalytic DC-weight sequence. If $c : \R,0 \to V \cq G$ is a germ of a $C^M$ curve and 
$\bar c : \R,0 \to \Si$ is a $C^\infty$ germ such that $\si \o \bar c = c$, then $\bar c$ is a $C^M$ germ.
\end{enumerate}
\end{theorem}

\begin{demo}{Proof}
Apply the results for finite group representations to the $W$-module $\Si$.
\qed\end{demo}

\begin{remark}
The conclusion in $(4)$ is wrong if $\bar c$ is a $C^\infty$ germ in $V$ which lies not in a Cartan subspace. 
See remark \ref{rmkp} and \ref{rmk1}. 
\end{remark}


\def\cprime{$'$}
\providecommand{\bysame}{\leavevmode\hbox to3em{\hrulefill}\thinspace}
\providecommand{\MR}{\relax\ifhmode\unskip\space\fi MR }
\providecommand{\MRhref}[2]{%
  \href{http://www.ams.org/mathscinet-getitem?mr=#1}{#2}
}
\providecommand{\href}[2]{#2}

\end{document}